\documentclass{amsart}
\usepackage{amssymb,amscd,verbatim,latexsym}
\usepackage{enumerate}
\usepackage[mathscr]{eucal}

\usepackage{color}
\newcommand\ede{ \, := \, }

\newcommand\rp{'}

\newcommand\dlp{``}
\newcommand\drp{''}

\newcommand\pullback{\sp{\downarrow\downarrow}}
\newcommand{\tto}{\rightrightarrows}
\newcommand\mathbfPsi{\mathbf \Psi}

\newcommand{\CC}{\mathbb C}

\newcommand{\RR}{\mathbb R}

\newcommand{\maB}{\mathcal B}
\newcommand{\maC}{\mathcal C}

\newcommand{\maF}{\mathcal F}
\newcommand{\maG}{\mathcal G}
\newcommand{\maH}{\mathcal H}

\newcommand{\maK}{\mathcal K}
\newcommand{\maL}{\mathcal L}

\newtheorem{theorem}{Theorem}[section]
\newtheorem{proposition}[theorem]{Proposition}

\theoremstyle{definition}
\newtheorem{definition}[theorem]{Definition}
\theoremstyle{remark}

\newtheorem{example}[theorem]{Example}

\author[V. Nistor]{Victor Nistor} \address{Universit\'{e} de Lorraine,
  UFR MIM, Ile du Saulcy, CS 50128, 57045 METZ, France and
  Inst. Math. Romanian Acad.  PO BOX 1-764, 014700 Bucharest Romania}
\email{victor.nistor@univ-lorraine.fr}

\date\today
\thanks{V.N. has been partially supported by
  ANR-14-CE25-0012-01.\\
Manuscripts available from {\bf http:{\scriptsize
    //}iecl.univ-lorraine.fr{\scriptsize
    /}$\tilde{}$Victor.Nistor{\scriptsize /}}}

\begin{document}

\title[Groupoids]{Fredholm criteria for pseudodifferential operators
and induced representations of groupoid algebras}

\begin{abstract}  
We characterize the groupoids for which an operator is Fredholm if,
and only if, its principal symbol and all its boundary restrictions
are invertible. A groupoid with this property is called {\em
  Fredholm}. Using results on the Effros-Hahn conjecture, we show that
an almost amenable, Hausdorff, second countable groupoid is Fredholm.
Many groupoids, and hence many pseudodifferential operators appearing
in practice, fit into this framework. We show that the
desingularization of groupoids preserves the class of Fredholm
groupoids.
\end{abstract}

\maketitle
\tableofcontents

\section{Introduction}

We obtain necessary and sufficient conditions for operators modeled
by groupoids to be Fredholm. Examples include operators obtained by
desingularization of singular spaces by successively blowing up the
lowest dimensional singular strata. We begin with a general study of
Fredholm conditions for pseudodifferential operators in the framework
of Fredholm groupoids.  A {\em Fredholm groupoid} is, by definition, a
locally compact groupoid with a Haar system for which the Fredholm
property is equivalent to the invertibility of the principal symbol
and of its fiberwise boundary restrictions. We obtain a general
characterization of Fredholm groupoids.  In particular, using some
results of Renault \cite{renault87, renault91} and Ionescu and
Williams \cite{ionescuWilliamsEHC}, we show that an almost amenable,
second-countable, Hausdorff groupoid is Fredholm.

Let $\maG$ be a groupoid with base $M$ modeling the analysis on some
singular space.  An {\em $A(\maG)$-tame submanifold $L \subset M$} is
one that has, by definition, a tubular neighborhood on which $A(\maG($
becomes a pull-back Lie algebroid.  The \dlp
desingularization\drp\ $[[\maG:L]]$ of $\maG$ along $L$
\cite{nistorDesing} is the a groupoid model-ling the analysis on the
space obtained by blowing-up $L$.  The space of units of the
desingularization $[[\maG:L]]$ is $[M:L]$, the blow-up of $M$ along
$L$. The desingularization groupoid is not a blown-up space, however.
We use the explicit structure of the desingularized groupoid
$[[\maG:L]]$ (see \cite{nistorDesing}) to show that it is Fredholm if
$\maG$ is.  Our results specialize to yield Fredholm conditions for
operators on manifolds with cylindrical and poly-cylindrical ends, on
manifolds that are asymptotically Euclidean or asymptotically
hyperbolic, on products of such manifolds, on manifolds that locally
at infinity are products of such manifolds, and on others. Most of the
(generally) easy proofs are contained in \cite{nistorFredCondLong},
this paper being a summary and update of some results in that paper.
We thank Ingrid and Daniel Belti\c t\u a, Claire Debord, Siegfried
Echterhoff, Vladimir Georgescu, Marius M\u antoiu, Jean Renault,
Steffen Roch, Georges Skandalis, and Dana Williams for useful
discussions.

\section{Fredholm groupoids}
Recall that a {\em groupoid} $\maG$ is a small category in which every
morphism is invertible. We shall write $\maG \tto M$ for a groupoid
with objects (or {\em units}) $M$. The domain and range of a morphism
therefore give rise to maps $d, r : \maG \to M$.  We refer to
\cite{MoerdijkMrcun02, renaultBook} for the results and concepts
used--but not recalled--in this paper.
Let $\maG \tto M$ be a locally compact groupoid endowed with a Haar
system $(\lambda_x)_{x \in M}$. We denote by $C\sp{\ast}(\maG)$ the
$C\sp{\ast}$-algebra of $\maG$ and by $C\sp{\ast}_r(\maG)$ the {\em
  reduced} $C\sp{\ast}$-algebra of $\maG$.  Also, we denote $\maG_A :=
d\sp{-1}(A)$ and $\maG_A\sp{B} := d\sp{-1}(A) \cap r\sp{-1}(B)$.  If
$A$ is invariant, in the sense that $\maG_{A} = \maG\sp{A} :=
r\sp{-1}(A)$, then $\maG_A$ is a groupoid. As usual, we associate to
any $x\in M$ the \emph{regular} representation $\pi_x\, \colon\,
C\sp{\ast}(\maG) \to \maB(L^2(\maG_x,\lambda_x))$
\begin{equation*}
  \pi_x(\varphi)\psi(g) \ede \varphi * \psi(g) \ede \int_{\maG_{d(g)}}
  \varphi(gh^{-1}) \psi(h) d\lambda_{d(g)}(h) \,, \quad \phi \in
  \maC_c(\maG) \,.
\end{equation*} 
All the morphisms (and representations) in this paper will preserve
the involution.

\begin{example}\label{ex.pullback}
Recall that the {\em pair groupoid} $\maH := A \times A$ is the
groupoid having {\em exactly} one arrow between any two units. Let
$\maG \tto L$ be a groupoid and $f : M \to L$. An important
generalization of the pair groupoid is the {\em fibered pull-back
  groupoid}:
\begin{equation*}
  f\pullback (\maG) \ede \{\, (m, g, m\rp) \in M \times \maG \times M,
  \ f(m) = r(g),\, d(g) = f(m\rp) \, \} \,,
\end{equation*}
with units $M$ and product $(m, g, m\rp) (m\rp, g\rp, m\drp) = (m, g
g\rp, m\drp)$.
\end{example}

\begin{definition}\label{def.fredholm}
A locally compact groupoid $\maG \tto M$ is called {\em Fredholm} if:
\begin{enumerate}[(i)]
\item $\maG_{M_0}\simeq M_0\times M_0$ for some open, dense,
  $\maG$-invariant subset $M_0\subset M$.
\item For any $a \in C^{\ast}_r(\maG)$, we have that $1 +
  \pi_{x_0}(a)$, $x_0 \in M_0$, is Fredholm if, and only if, $1 +
  \pi_x(a)$ is invertible for any $x \in F:= M \smallsetminus M_0$.
\end{enumerate}
\end{definition}
  
Both $M_0$ and $F$ are uniquely determined by $\maG$, so this notation
will remain fixed in what follows. Also, in Definition
\ref{def.fredholm}, all representations $\pi_{x_0}$, $x_0 \in M_0$,
are unitarily equivalent to the {\em vector representation} $\pi_0 :
C\sp{\ast}(\maG) \to \maL(L\sp{2}(M_0))$ obtained by identifying $r:
\maG_{x_0} \simeq M_0$.

Recall \cite{Roch} that if $A$ is a $C\sp{\ast}$-algebra with unit,
then a set $\maF$ of representations of $A$ is called {\em
  invertibility sufficient} if the following condition is satisfied:
\dlp $a \in A$ is invertible if, and only if, $\phi(a)$ is invertible
for all $\phi \in \maF$.\drp\ If $A$ does not have a unit, we replace
$A$ with $A\sp{+} := A \oplus \CC$ and $\maF$ with $\maF\sp{+} := \maF
\cup \{\chi_0 : A\sp{+} \to \CC\}$ \cite{nistorPrudhon}.  The
following two results give a first characterization of Fredholm
groupoids.

\begin{theorem}  \label{thm.Fredholm.Cond}
Let $\maG \tto M$ be a Fredholm groupoid. Then
\begin{enumerate}[(i)]
\item The vector representation $\pi_0 \colon C_r\sp{\ast}(\maG) \to
  \maL(L^2(M_0))$ is injective.
\item The canonical projection induces an isomorphism
  $C_r\sp{\ast}(\maG)/C_r\sp{\ast}(\maG_{M_0}) \simeq
  C_r\sp{\ast}(\maG_F)$.
\item $\{\pi_x,\, x \in F\}$ is an invertibility sufficient set of
  representations of $C_r\sp{\ast}(\maG_F)$.
\end{enumerate}
\end{theorem}

For index theory, the first two conditions of the theorem are
enough. Invertibility sufficient families of representations consist
of non-degenerate representations. A non-degenerate representation of
a (closed, two-sided) ideal in a $C\sp{\ast}$-algebra has a unique
extension to the whole algebra.

The following strong converse of Theorem \ref{thm.Fredholm.Cond} holds
true.

\begin{theorem} \label{thm.Fredhol.Cond.Conv}
Let $\maG \tto M$ be a locally compact groupoid satisfying the three
conditions (i-iii) of Theorem \ref{thm.Fredholm.Cond}. Then, for any
unital $C^{\ast}$-algebra $\mathbfPsi$ containing $C^{\ast}_r(\maG)$
as an essential ideal and for any $a \in \mathbfPsi$, we have that
$\pi_0(a)$ if Fredholm if, and only if, $\pi_x(a)$ is invertible for
each $x \notin M_0$ {\bf and} the image of $a$ in
$\mathbfPsi/C^{\ast}_r(\maG)$ is invertible.
\end{theorem}

\section{Relation to the Effros-Hahn conjecture} 
We now want to obtain some more concrete and easier to use conditions
for a groupoid $\maG$ to be Fredholm. We shall say that a locally
compact groupoid $\maG$ {\em has the weak inclusion property}
(wi-property, for short) if every irreducible representation of
$C\sp{\ast}(\maG)$ is {\em weakly contained} in a representation of
$C\sp{\ast}(\maG)$ induced from a representation of an isotropy
subgroup $\maG_y\sp{y} := d\sp{-1}(y) \cap r\sp{-1}(y)$
\cite{ionescuWilliamsEHC, renault91} (equivalently, if every primitive
ideal of $C\sp{\ast}(\maG)$ {\em contains} an ideal induced from a
representation of an isotropy subgroup $\maG_y\sp{y}$).  A groupoid
$\maG$ with the wi-property and such that all the groups $\maG_y^y$
are amenable will be called EH-{\em amenable}. Recall that a locally
compact groupoid $\maG$ {\em satisfies the generalized Effros-Hahn
  (EH) conjecture} if every primitive ideal of $C\sp{\ast}(\maG)$ is
{\em induced} from a representation of an isotropy group
$\maG_y\sp{y}$ \cite{EchterhoffMemoirs, ionescuWilliamsEHC,
  ionescuWilliamsIndRep, renaultBook, williamsBook}.

\begin{example}\label{ex.EH-amenable}
 Let $\maG \tto B$ be a locally trivial bundle of groups (so $d = r$)
 with typical fiber a locally compact group $G$.  Also, let $f : M \to
 B$ be a continuous map that is a local fibration. Then $f\pullback
 (\maG)$ is a locally compact groupoid with a Haar system that
 satisfies the generalized EH conjecture, and hence it has the
 wi-property. It will be EH-amenable if, and only if, the group $G$ is
 amenable.
\end{example}

Recall that a groupoid is called {\em metrically amenable} if the
canonical surjection $C\sp{\ast}(\maG) \to C_r\sp{\ast}(\maG)$ is
injective \cite{SimsWilliams}.  We shall need two results from
\cite{nistorPrudhon} (see also \cite{ExelInvGr}).

\begin{proposition}\label{prop.sufficient}
Let $\maG\tto F$ be an $EH$-amenable locally compact groupoid.  Then
the family of regular representations $\{\pi_y, y \in F\}$ of
$C\sp{\ast}(\maG)$ is invertibility sufficient. In particular, $\maG$
is metrically amenable.
\end{proposition}

The class of EH-amenable groupoids is closed under extensions.

\begin{proposition}\label{prop.EHexact}
Let $\maG \tto M$ be a locally compact groupoid, $M_0 \subset M$ be a
$\maG$-invariant open subset, and $F := M \smallsetminus M_0$. Then
$\maG$ is EH-amenable if, and only if, both $\maG_F$ and $\maG_{M_0}$
are EH-amenable. The same holds if one replaces \dlp is
EH-amenable\drp\ with \dlp satisfies the generalized EH
conjecture\drp\ or \dlp has the wi-property.\drp
\end{proposition}

The following result leads to more applicable Fredholm conditions.

\begin{proposition}\label{prop.amenable}
 Let $\maG \tto M$ be a Hausdorff, locally compact groupoid with an
 open, dense, $\maG$-invariant subset $M_0 \subset M$ such that
 $\maG_{M_0} \simeq M_0 \times M_0$. If $\maG$ is EH-amenable, then
 $\maG$ is Fredholm.
\end{proposition}

Let $U_i \subset U_{i+1} \subset M$ be open, $\maG$-invariant subsets
of $M$, with $U_{-1} = \emptyset$ and $U_{N} = M$. If
$\maG_{U_{i+1}\smallsetminus U_i}$ is topologically amenable for all
$i$, then we shall say that the locally compact groupoid $\maG \tto M$
is {\em almost amenable}.  By combining the above two propositions
with the proof of the generalized EH conjecture
\cite{ionescuWilliamsEHC, renault87, renault91} for amenable,
Hausdorff, second countable groupoids, we obtain the following result.

\begin{theorem}\label{thm.amenable}
Let $\maG \tto M$ be an almost amenable, Hausdorff, second countable
groupoid.  Then $\maG$ satisfies the generalized EH conjecture. If
also $\maG_{M_0} \simeq M_0 \times M_0$ for an open, dense,
$\maG$-invariant subset $M_0 \subset M$, then $\maG$ is Fredholm.
\end{theorem}

We are interested in Fredholm groupoids because of their applications
to Fredholm conditions. Let $\maG$ be a continuous family groupoid
\cite{LMN1} and $\Psi\sp{m}(\maG)$ be the space of order $m$,
classical pseudodifferential operators $P = (P_x)_{x \in M}$ on $\maG$
\cite{LMN1} (see \cite{ASkandalis2, aln2, Monthubert, NWX,
  vanErpYuncken} for Lie groupoids, which are continuous family
groupoids).  Recall that, by definition, each $P_x \in
\Psi\sp{m}(\maG_x)$, $x \in M$.  Also, $P$ acts on $M_0$ via $P_{x_0}:
H^s(M_0) \to H^{s-m}(M_0)$, $x_0 \in M_0 = r(\maG_{x_0}) \simeq
\maG_{x_0}$. The following result is interesting only in the case $M$
compact.

\begin{theorem}\label{thm.nonclassical2}
Let $\maG \tto M$ be a Fredholm, continuous family groupoid and let
$M_0 \subset M$ be the dense, $\maG$-invariant subset such that
$\maG_{M_0} \simeq M_0 \times M_0$.
We have
\begin{multline*} 
	P : H^s(M_0) \to H^{s-m}(M_0), \ P \in \Psi\sp{m}(\maG)\,,
        \mbox{ is Fredholm} \ \ \Leftrightarrow \ \ P \mbox{ is
          elliptic and }\\
	\ P_{x} : H^s(\maG_x) \to H^{s-m}(\maG_x) \ \mbox{ is
          invertible for all } x \in F := M \smallsetminus M_0 \,.
\end{multline*}
\end{theorem}

This theorem is proved by considering $a := (1 + \Delta)\sp{(s-m)/2} P
(1 + \Delta)\sp{-s/2}$, which belongs to the closure of
$\Psi\sp{0}(\maG)$, by the results in \cite{LNGeometric}. For the next
theorem, however, one has to consider the Cayley transform of $P$
instead of the operator $a$.
 
\begin{theorem}\label{thm.nonclassical.sp}
Let $\maG \tto M$ be as in Theorem \ref{thm.nonclassical2} and let $P
\in \Psi\sp{m}(\maG)$ be an elliptic operator. Then its essential
spectrum is
\begin{equation*} 
	\sigma_{ess}(P) \ = \
	\begin{cases}
	\ \ \cup_{x \in F}\, \sigma(P_{x})\ \ \ & \mbox{ if } m > 0\\
	\ \cup_{x \in F}\, \sigma(P_{x}) \  \cup \ Im(\sigma_0(P)) & 
	\mbox{ if } m \le 0\,.
	\end{cases}
\end{equation*}
\end{theorem}
 
The above two theorems extend to operators acting on vector bundles on
$M$.  The operators $P_x$ are the analogues in our setting of the
``limit operators'' considered in \cite{DamakGeorgescu, RochBookLimit}
and many other references.  See also \cite{beltitzas15, DLR,GI2,
  GeorgescuNistor1, LMN1, LauterMoroianu1, LNGeometric, MantoiuReine, 
  Mantoiu1, MazzeoMelroseAsian, Schulze, Rab3, RochBookNGT, schulzeWong2014,
  SchroheFrechet, SchSch1} and the references therein for related
results.

\section{Desingularization and Fredholm conditions}

We want an ample supply of Fredholm groupoids. In this section, we
recall the desingularization procedure along a \dlp
tame\drp\ submanifold of the set of units of a Lie groupoid
\cite{nistorDesing}. Recall the thick pull-back $\pi\pullback (B)$ of
Example \ref{ex.pullback}.

\begin{definition}\label{def.tame}
 Let $A \to M$ be a Lie algebroid.  An {\em $A$-tame submanifold} $L
 \subset M$ is a submanifold that has a tubular neighborhood $L
 \subset U$, such that there exists a Lie algebroid $B \to L$ and an
 isomorphism $A \vert_{U} \simeq \pi\pullback(B)$ over $U$, $\pi : U
 \to L$.
\end{definition}

By $A(\maG)$ we denote the Lie algebroid of a Lie groupoid $\maG$. We
have the following structure result in a neighborhood of a tame
submanifold.

\begin{theorem}\label{thm.tame}
 Let $\maG \tto M$ be a Lie groupoid and let $L \subset M$ be an
 $A(\maG)$-tame submanifold of $M$ with tubular neighborhood $\pi : U
 \to L \subset U$.  Assume that the fibers of $\pi : U \to L$ are
 simply-connected. Then the reduction groupoids $\maG_{L}\sp{L}$ and
 $\maG_U\sp{U}$ are Lie groupoids and there exists a natural
 isomorphism over $U$ of Lie groupoids
\begin{equation*}
  \maG_{U}\sp{U} \, \simeq \, \pi\pullback ( \maG_{L}\sp{L} ) \ede U
  \times_{\pi} \maG_{L}\sp{L} \times_{\pi} U \,.
\end{equation*}
\end{theorem}

From now on, we shall assume that $\maG \tto M$ and $L \subset M$ are
as in Theorem \ref{thm.tame} and use this theorem to define the
desingularization $[[\maG:L]]$.  To construct $[[\maG:L]]$, we proceed
in four steps \cite{nistorDesing}, the first three of which deal with
the particular case $\maG = f\pullback (\maH)$, where $f : M = U \to
L$ is a vector bundle. We denote by $S$ the set of unit vectors in $U$
for some fixed metric.

\subsection*{Construction of the desingularization}
{\em Step one.}  We first consider the adiabatic groupoid $\maH_{ad}$
of $\maH$ \cite{ConnesBook, DebordSkandalis, NWX}. It is a Lie
groupoid with units $L \times [0, \infty)$ and Lie algebroid
  $A(\maH_{ad}) = A(\maH) \times [0, \infty) \to L \times [0,
      \infty)$, which, as a vector bundle, is the pull-back of
      $A(\maH) \to L$ to $L \times [0, \infty) \to L$.  The Lie
        algebroid structure on the sections of $A(\maH_{ad})$ is not
        that of a pull-back, but is given by $[X, Y](t) = t [X(t),
          Y(t)].$ As a set, $\maH_{ad}$ is the disjoint union
\begin{equation*}
  \maH_{ad} \, := \, A(\maH) \times \{0\} \, \sqcup \, \maH \times (0,
  \infty) \,.
\end{equation*}
The groupoid structure of $\maH_{ad}$ is such that $A(\maH) \times
\{0\}$ has the Lie groupoid structure of a bundle of Lie groups and
$\maG \times (0, \infty)$ has the product Lie groupoid structure (that
is $(0, \infty)$ has only units, and all orbits are reduced to a
single point).

{\em Step two.}\ Let $\pi : S \to L$ be the projection. We denote also
by $\pi$ the resulting map $S \times [0, \infty) \to L \times [0,
    \infty)$.  Then we consider the Lie groupoid
    $\pi\pullback(\maH_{ad})$.

{\em Step three.}\ Let $\RR_{+}\sp{*} = (0, \infty)$ act by dilations
on the $[0, \infty)$ variable on $\pi\pullback (\maH_{ad})$ and
  consider the semi-direct product $\pi\pullback (\maH_{ad})\rtimes
  \RR_{+}\sp{*}$ \cite{DebordSkandalis}. We then define
\begin{equation}\label{eq.first}
 [[\maG:L]] := \pi\pullback (\maH_{ad})\rtimes \RR_{+}\sp{*}\,,
 \ \mbox{ if } \maG = \pi \pullback(\maH) \mbox{ for a vector bundle }
 \pi : M \to L \,.
\end{equation}
{\em Step three.}\ Let $\maG \tto M$ be a Lie algebroid and $L \subset
M$ be an $A(\maG)$-tame submanifold. Let $W := M \smallsetminus L$ and
let $U$ be as in Theorem \ref{thm.tame}.  Then we glue $\maG_W\sp{W}$
and $[[\maG_U\sp{U} : L]]$ along the common open subset
\cite{Gualtieri, nistorDesing} to obtain $[[\maG:L]]$.
  
If $\maG$ is Hausdorff, then $[[\maG_U\sp{U} : L]]$ is also Hausdorff,
\cite{nistorDesing}. We denote by $[M:L]$ the {\em blow-up} of $M$
with respect to $L$, it is obtained by replacing $L$ with the set of
unit vectors $S$ of its normal bundle in $M$. So $S = [M:L]
\smallsetminus (M \smallsetminus L) = [M:L] \smallsetminus W$.  Let
$\maG$ be as in Theorem \ref{thm.tame}.  We then have the following
structural result for the desingularization~$\maK := [[\maG:L]]$.

\begin{proposition}\label{prop.blow-up} 
The desingularization $\maK := [[\maG:L]]$ is a Lie groupoid with
units $[M: L]$. The subset $S \subset [M:L]$ is closed and
$\maK$-invariant. The restriction $\maK_{S} = [[\maG: L]]
\smallsetminus \maG_W\sp{W}$ is isomorphic to the fibered pull-back
$\pi\pullback (A(\maK) \rtimes \RR_{+}\sp{*})$ to $S$ via the natural
projection $\pi : S \to L$, where $A(\maK) \rtimes \RR_{+}\sp{*}$ is
regarded as a bundle of Lie groups.  The inclusion $\maG_W\sp{W} \to
\maG$ induces an isomorphism $C\sp{\ast}(\maG_W\sp{W}) \simeq
C\sp{\ast}(\maG)$.
\end{proposition}

One sees that the resulting glued set is a Hausdorff groupoid as in
\cite{Gualtieri}.

\begin{theorem}\label{thm.Fred.desing}
 Let $\maG \tto M$ be a Lie groupoid. Let us assume that $\maG$ is
 Fredholm.  Let $L \subset M$ be an $A(\maG)$-tame submanifold. Then
 $[[\maG:L]]$ is also Fredholm.
\end{theorem}

Desingularization preserves the class of Fredholm groupoids.

\begin{theorem}\label{thm.final}
Let us assume that $\maG$ is obtained from a pair groupoid $M \times
M$ (with $M$ smooth) by a sequence of desingularizations with respect
to tame submanifolds. Then $\maG$ is a Hausdorff Fredholm Lie
groupoid.
\end{theorem}

This theorem can be used to obtain Fredholm conditions for operators
on polyhedral domains, as well as on some other stratified spaces.

\def\cprime{$'$}

\end{document}